
\def\ifudf#1{\expandafter\ifx\csname #1\endcsname\relax}
\newif\ifpdf \ifudf{pdfoutput}\pdffalse\else\pdftrue\fi
\ifpdf\pdfpagewidth=210mm \pdfpageheight=297mm
 \else\fi
\hsize=149mm \vsize=245mm \hoffset=5mm \voffset=0mm
\parskip=1ex minus .3ex \parindent=0pt
\everydisplay={\textstyle}
\font\tbb=bbmsl10 \font\sbb=bbmsl10 scaled 700
 \font\fbb=bbmsl10 scaled 500
\newfam\bbfam \textfont\bbfam=\tbb
 \scriptfont\bbfam=\sbb \scriptscriptfont\bbfam=\fbb
\def\bb{\fam\bbfam}%
\font\tfk=eufm10 \font\sfk=eufm7 \font\ffk=eufm5
\newfam\fkfam \textfont\fkfam=\tfk
 \scriptfont\fkfam=\sfk \scriptscriptfont\fkfam=\ffk
\def\fk{\fam\fkfam}%
\def\fn[#1]{\font\TmpFnt=#1\relax\TmpFnt\ignorespaces}
\def\em{\expandafter\ifx\the\font\tensl\rm\else\sl\fi}
\def\dt{\number\day.\number\month.\number\year/\the\time}
\def\\{\hfill\break}
\newtoks\secNo \secNo={0}
\newcount\ssecNo \ssecNo=0
\def\section #1#2\par{\goodbreak\vskip 3ex\noindent
 \global\secNo={#1}%
 \global\ssecNo=0 \global\eqnNo=0
 {\fn[cmbx10 scaled 1200]#1#2}\vglue 1ex}
\def\subsection #1.{\goodbreak\vskip 3ex\noindent
 \global\advance\ssecNo by 1
 {\fn[cmbx10]\the\secNo.\the\ssecNo.~#1}\vglue 1ex}
\def\url #1<#2>{\ifpdf\pdfstartlink
 attr {/Border [9 9 1] /C [1 0.7 0.7]}%
 user {/Subtype /Link /A << /S /URI /URI (#2) >>}\fi%
 {#1}\ifpdf\pdfendlink\fi\ignorespaces}
\def\msc#1{\url #1 <https://zbmath.org/classification/?q=cc:#1>}

\def\href#1<#2>{\leavevmode
 \ifpdf\pdfstartlink attr {/Border [0 0 0 ]} goto name {#2}\fi
 {#1}\ifpdf\pdfendlink\fi\ignorespaces}
\def\label@#1:#2@{\ifudf{#1}
 \expandafter\xdef\csname#1\endcsname{#2}\else
 \errmessage{label #1 already in use!}\fi}
\label@liemob:1.1@
\label@kilsurf:1.2@
\label@rodrigues:1.3@
\label@cayham:1.4@
\label@kilhgm:1.5@
\label@eucisom:1.6@
\label@eucsph:1.7@
\label@poisurf:1.8@
\label@poihgm:1.9@
\label@poihsc:1.10@
\label@rconfak:Lemma 1@
\label@bryantlw:2.1@
\label@isocong:2.2@
\label@lwcharpf:2.3@
\label@lwchar:Lemma 2@
\label@mobsc:2.4@
\label@paralhf:2.5@
\label@iscong:2.6@
\label@sgauss:2.7@
\label@curvlw:Thm 3@
\label@isocoord:3.1@
\label@radiusbc:3.2@
\label@bcm:Thm 5@
\label@fig:Fig 1@
\label@2026bjkp:1@
\label@1923bi:2@
\label@1929bl:3@
\label@1997bjpp:4@
\label@2010bjr:5@
\label@2012bjr:6@
\label@2019bjlm:7@
\label@1900ca:8@
\label@2008ce:9@
\label@2025cnopry:10@
\label@2004gmm:11@
\label@2001jmn:12@
\label@2003imdg:13@
\label@2009ko:14@
\label@2002liro:15@
\label@2020pe3:16@
\label@2020pe:17@
\label@1994sm:18@
\label@2020ze:19@
\def\@#1:#2@{\ifpdf\pdfdest name {#1} xyz\fi {#2}}
\def\:#1:{\href\ifudf{#1}??\else\csname#1\endcsname\fi<#1>}
\newcount\refNo \refNo=0
\def\refitem#1 {\global\advance\refNo by 1
 \item{\@#1:\number\refNo@.}}
\let\pleqno\eqno \newcount\eqnNo \eqnNo=0
\def\eqno#1$${\global\advance\eqnNo by 1
 \pleqno{\rm(\@#1:\the\secNo.\number\eqnNo@)}$$}
\def\R{{\bb R}} 
\def\C{{\bb C}} 
\def\P{{\bb P}}
\def\lc{{\cal L}}
 
\def\del#1{{\partial\over\partial #1}}

\newtoks\title \newtoks\stitle \newtoks\author \newtoks\sauthor
\title={A Bianchi-Cal\`o method for Bryant type surfaces}
\stitle={Bianchi-Cal\`o method}
\author={F Burstall, U Hertrich-Jeromin, G Szewieczek}
\sauthor={F Burstall et al}
\ifpdf\pdfoutput=1\pdfadjustspacing=1\pdfinfo{%
 /Title (\the\title) /Author (\the\author) /Date (\dt)}\fi
\headline={\ifnum\pageno>1
 {\fn[cmr7]\the\sauthor\hfill\the\stitle}\fi}
\footline={{\fn[cmr7]\hfil-\folio-\hfil}}
\centerline{{\fn[cmbx10 scaled 1440]\the\title}}\vglue .2ex
\centerline{{\fn[cmr7]\the\author}}\vglue 3em plus 3ex
\centerline{\vtop{\hsize=.8\hsize{\bf Abstract.}\enspace
 We present a Bianchi-Cal\`o type construction method for
 Bryant type linear Weingarten surfaces in hyperbolic space.
}}\vglue 2em
\centerline{\vtop{\hsize=.8\hsize{\bf MSC 2020.}\enspace
 {\it\msc{53C42}\/}, {\it\msc{53A10}\/},
  \msc{53A35}, \msc{37K25}, \msc{37K35}
}}\vglue 1em
\centerline{\vtop{\hsize=.8\hsize{\bf Keywords.}\enspace
 linear Weingarten surface; Bryant type surface;
 isothermic surface; sphere congruence;
 hyperbolic geometry; Euclidean geometry; M\"obius geometry;
 Lie sphere geometry;
 Weierstrass type representation.
}}\vglue 3em

\section Introduction

Bianchi [\:1923bi:, \S496(b) (p.612)]
 provides a method to construct a cmc-$1$ surface in hyperbolic
 geometry from a given holomorphic hyperbolic Gauss map
 ${\fk h}:\C\supset{\cal D}\to\C$,
 introduced
 in [\:1923bi:, \S495 (p.607)];
in fact, he provides an explicit and integration free way
 to construct a particular Darboux transform of this
 hyperbolic Gauss map, cf [\:2001jmn:, Sect 4.1],
by giving an explicit formula for the (Euclidean) centre
 surface of the enveloped horosphere congruence
 in Poincar\'e half space,
 in [\:1923bi:, \S328 (p.133)]:
$$
  {\fk c} = (r,{\fk h}):\C\supset{\cal D}\to\R\times\C
   \enspace{\rm with}\enspace
  r = {(1+|z|^2)\,|{\fk h'}|\over 2}.
$$
This representation of the centre surface is attributed
 to Cal\`o [\:1900ca:],
 as the solution to a rolling problem,
that is,
 a problem of determining isometric surface pairs.
Thus we follow the modern account [\:2002liro:] of this method
 in calling it the ``Bianchi-Cal\`o method''.

The classical Bianchi-Cal\`o method features an integration
 free representation for cmc-$1$ surfaces in hyperbolic space,
similar to Small's representation [\:1994sm:\ (3.1)]
 in the Poincar\'e ball model
that the authors relate to their explicit parametrization
 [\:2002liro:\ (1)] obtained by the Bianchi-Cal\`o method
 in [\:2002liro:, Sect 3].

A second remarkable attribute of the Bianchi-Cal\`o method
 is the interplay of geometries involved, cf [\:2002liro:]:
this is not just an entertaining feature,
but it is an underlying cause for the workings of the
 Bianchi-Cal\`o method.
As the (Euclidean) centre surface of a (hyperbolic) horosphere
 congruence plays a role,
a relation between these two geometries is obviously paramount.
However, we shall demonstrate that also relations with
 sphere geometries play a key role
 in the (generalized) Bianchi-Cal\`o method,
as the appearance of the Darboux transformation
 in [\:1923bi:, \S496] suggests, cf [\:2020ze:]:
this analysis is the main focus of the present text.

The scope of representations in terms of holomorphic data
 for surfaces in hyperbolic space has in [\:2004gmm:] been
 extended to surfaces satisfying a relation (\:bryantlw:)
 of their Gauss and mean curvatures,
the ``Bryant type linear Weingarten surfaces'',
see also [\:2009ko:, Thm 3.8] and [\:2020pe:, Sect 4.4].
We observed in [\:2012bjr:, Sect 4] that these surfaces are
 characterized by the existence of an enveloped isothermic
 horosphere congruence,
which gives rise to a holomorphic hyperbolic Gauss map.
Thus the observations
 of [\:2001jmn:, Sect 4.1] and
 of [\:2009ko:, Thm 3.8]
suggested that a construction of Bryant type linear Weingarten
 surfaces via special Ribaucour transformation of
 a prescribed holomorphic hyperbolic Gauss map
 should be possible.

Indeed, in this text we observe (\:curvlw:) that
 said horosphere congruence envelops a Bryant type surface
 with parameter $\mu$ precisely when its induced metric
 has curvature $K_{(ds,ds)}=-\mu$,
and (\:rconfak:) that
 its radius function occurs as a conformal factor;
these ingredients then yield
our generalization (\:bcm:) of the Bianchi-Cal\`o method,
 where the horosphere congruence has centre surface
$$
  {\fk c} = (r,{\fk h}):\C\supset{\cal D}\to\R\times\C
   \enspace{\rm with}\enspace
  r = {(1-\mu|z|^2)\,|{\fk h'}|\over 2}.
$$

The text is organized as follows:
Sect 1 sets the scene, in particular, it fixes notions and
 constructions of various geometries that make an
 appearance;
Sect 2 then recalls relevant results from [\:2012bjr:] and
 establishes a relation between the Bryant type
 condition and the curvature of an enveloped horosphere
 congruence;
Sect 3 is devoted to our generalization of the Bianchi-Cal\`o
 method.

{\it Acknowledgements.\/}
This paper has taken 25 years to get written;
consequently, numerous discussions and ideas have influenced
 the work and
it seems impossible to mention all colleagues involved in
 those discussions.
However, we wish to mention some milestones of
 this 25-year endeavour:
{\parindent=2em
\item{$\bullet$} the second author got first interested in the
 Bianchi-Cal\`o method through [\:2002liro:] in the context of
 the work [\:2001jmn:] on cmc-$1$ surfaces in hyperbolic space
 --- here we would like to thank
  {\em E~Musso\/},
  {\em L~Nicolodi\/} and
  {\em W~Rossman\/}
 for discussions on cmc-$1$ surfaces in hyperbolic space,
 as well as
  {\em L~de~Lima\/} and {\em P~Roitman\/}
 for their inspiring paper;
\item{$\bullet$} the first and second author later got
 interested in flat fronts and more general Bryant type
 linear Weingarten surfaces in the context of Lie sphere
 geometry, see [\:2010bjr:, \:2012bjr:]
 --- and we would like to thank
  {\em S~Fujimori\/},
  {\em M~Kokubu\/},
  {\em M~Pember\/},
  {\em W~Rossman\/},
  {\em M~Umehara\/} and
  {\em K~Yamada\/}
 for inspiration and information through many enjoyable
 discussions;
\item{$\bullet$} the second and third author took the topic
 up more recently in the context of the MSc thesis [\:2020ze:]
 of
  {\em J~Zeilinger\/}
 whom we would like to thank for her interest in and
 our discussions about the Bianchi-Cal\`o method and
 Darboux transformations;
\item{$\bullet$} finally, joining forces, the present authors
 of this paper arrived at a reasonable understanding of
 the method that allows for generalization
 --- with thanks for further input and discussions about
 the topic to
  {\em J~Cho\/}
 and several of the aforementioned colleagues.
\par}

We also gratefully acknowledge financial support from
the Austrian Science Fund FWF through the FWF/JSPS joint
research project I3809 ``Geometric shape generation''.
The figures were created using the computer algebra system
{\tt Mathematica}.

\section 1. Setting the scene

We are interested in the relations between the geometry
 of a surface (in hyperbolic space),
 of its (hyperbolic) Gauss map and
 of the enveloped (horo-)sphere congruence enveloped by both,
 given in terms of its (Euclidean) centre surface.
For the investigation of such relations,
it will be convenient to alternate between different viewpoints
and ambient geometries:
 hyperbolic geometry,
 M\"obius geometry as well as
 Lie sphere geometry.
For the convenience of the reader,
we shall start with a concise discussion of those viewpoints;
for more details the reader is referred to
[\:1929bl:; \:2003imdg:; \:2008ce:; \:2025cnopry:].

\subsection Sphere geometries.
We consider all involved geometries as subgeometries of
 {\em Lie sphere geometry\/}:
thus points of the {\em Lie quadric\/}
 $\P(\lc^5)\subset\P(R^{4,2})$
 are interpreted as $2$-spheres,
and a surface is described by a line congruence
 in the Lie quadric, that is,
 by a {\em Legendre map\/} into the Grassmannian
 of null $2$-planes in $\R^{4,2}$,
$$
  \langle\varphi,\tau\rangle:\Sigma^2\to{\rm Gr}_0(\R^{4,2}),
$$
 where $(\varphi,\varphi)=(\varphi,\tau)=(\tau,\tau)=0$
 and the {\em contact condition\/} $(d\varphi,\tau)=0$ holds.

Fixing a {\em point sphere complex\/} $p\in\R^{4,2}$
 with $(p,p)=-1$
yields a splitting
 $\R^{4,2}=\langle p\rangle\oplus_\perp\R^{4,1}$;
spheres $\langle\xi\rangle\in\P(\lc^4)\subset\P(\R^{4,1})$
 are now considered as points while
{\em proper spheres\/} $\langle\sigma\rangle\not\in\P(\lc^4)$
 can be described by a normalized lift $\sigma=p+s\in\lc^5$
 with $s\in S^{3,1}\subset\R^{4,1}$.
Thus a surface in {\em M\"obius geometry\/} is described by
 a point map $\langle f\rangle$ and
 an {\em enveloped sphere congruence\/} $t$,
$$
  f:=\varphi:\Sigma^2\to\lc^4
   \enspace{\rm and}\enspace
  t:=\tau-p:\Sigma^2\to S^{3,1},
\eqno liemob$$
where $\varphi\perp p$ and
 $t$ can be interpreted as a unit normal field of $f$,
 as $(df,t)=(d\varphi,\tau)=0$.

\subsection Hyperbolic geometry: Killing model.
Now fixing a {\em space form vector\/}
 $q\in\R^{4,1}=\{p\}^\perp\subset\R^{4,2}$
 with $(q,q)=1$
yields $\R^{4,2}=\langle p,q\rangle\oplus_\perp\R^{3,1}$;
the projective light cone $\P(\lc^3)\subset\P(\R^{3,1})$
 models a (M\"obius geometric) $2$-sphere,
that we consider as the ideal boundary of
 Killing's familiar hyperboloid model of hyperbolic $3$-space:
any {\em proper point\/} $\langle x\rangle\not\in\P(\lc^3)$ can
 be described by a normalized lift $x\in\lc^4$ with $(x,q)=-1$,
 that is,
 up to a shift $q+x\in\R^{3,1}$ proper points are in the
 standard $2$-sheeted hyperboloid of $\R^{3,1}$.
Furthermore, spheres $t\in S^{2,1}\subset\R^{3,1}$ intersect
 the ideal boundary sphere $q\in S^{3,1}$ orthogonally
hence are planes of the hyperbolic geometry
 given by $p,q\in\R^{4,2}$.
The geometry of a surface in hyperbolic geometry
 is now described in terms of
 its {\em hyperbolic lift\/} and
 its de Sitter space valued Gauss map,
 interpreted as its {\em tangent plane map\/},
$$
  f:\Sigma^2\to K^3:=\{x\in\lc^4\,|\,(x,q)=-1\}
   \enspace{\rm and}\enspace
  t:\Sigma^2\to\{s\in S^{3,1}\,|\,(s,q)=0\}.
\eqno kilsurf$$
As usual,
 curvature line coordinates $(u,v)$ and principal curvatures
 $k_1,k_2$ of $f$ can be characterized by Rodrigues' equations
$$
  0 = t_u+k_1f_u = t_v+k_2f_v,
\eqno rodrigues$$
and the mean and (extrinsic) Gauss curvatures of $f$ are given
 by the Cayley-Hamilton formula
$$
  0 = (dt,dt) + 2H\,(dt,df) + K\,(df,df).
\eqno cayham$$
Furthermore, 
the two hyperbolic Gauss maps
 $h^\pm=q+f\pm t:\Sigma^2\to\lc^3$
 of $f$ are
the second envelopes of two congruences
 $t\pm f:\Sigma^2\to S^{3,1}$
 of horospheres that envelop $f$,
  as $f,df\perp t\pm f$,
 and touch $\mp q$,
  as $(t\pm f,\mp q)=1$ resp $(\tau\pm\varphi,p\mp q)=0$.
We shall focus on one of these,
$$
  h := q+f+t:\Sigma^2\to\lc^3
   \enspace{\rm and}\enspace
  s := t+f = h-q:\Sigma^2\to S^{3,1}.
\eqno kilhgm$$
Note that $f$ and $h$ differentiate into the first derived
 bundle $s^{(1)}=\langle s\rangle\oplus ds(T\Sigma^2)$,
 hence are parallel normal sections of $s$,
 confirming that $s$ has flat normal bundle and, therefore,
 is a {\em Ribaucour sphere congruence\/};
more specifically,
 $h$ yields {\em Ribaucour coordinates\/} of $f$,
cf [\:2019bjlm:].

\subsection Euclidean geometry.
Similar to the way hyperbolic geometry arises as a subgeometry
 of M\"obius geometry by way of a light cone section,
so does Euclidean geometry:
consider
 $\R^{4,1}=\langle o,\infty\rangle\oplus_\perp\R^3$,
where $(o,\infty)$ denotes a pseudo-orthonormal basis 
 of the orthogonal space $\R^{1,1}=(\R^3)^\perp\subset\R^{4,1}$
 of the Euclidean $3$-space,
 that is, $(o,o)=(\infty,\infty)=0$ and $(o,\infty)=-1$;
then
$$
  \R^3\ni{\fk y}\to y :
   = o + {\fk y} + {({\fk y,y})\over 2}\,\infty
     \in\lc^4\subset\R^{4,1}
\eqno eucisom$$
yields an isometry of $\R^3$ onto a parabolic light cone
 section.
A sphere $s\in S^{3,1}$ that does not contain the point
 $\langle\infty\rangle\in\P(\lc^4)$,
 that is, $(s,\infty)\neq 0$,
can be written as
$$
  s = {1\over r}\,(o+{\fk c}+{({\fk c,c})-r^2\over 2}\,\infty)
      \in S^{3,1}
\eqno eucsph$$
in terms of its centre ${\fk c}\in\R^3$ and its radius $r\neq 0$;
note that $(s,y)=-{({\fk y-c,y-c})-r^2\over 2r}$ confirming that
 the point $y$ and the sphere $s$ are incident iff $(s,y)=0$.

\subsection Hyperbolic geometry: Poincar\'e half space.
We now combine hyperbolic and Euclidean geometries in
 our M\"obius geometric setting:
choosing $o,\infty$ and $q$ as before with $q\perp o,\infty$
 yields a decomposition
$\R^{4,1}=\langle o,q,\infty\rangle\oplus_\perp\R^2$,
 where $\R^2\cong\C$;
note that $q\in S^{3,1}$ represents the Euclidean plane
 $\C\subset\langle q\rangle\oplus\C\cong\R^3$
 (with a choice of orientation),
 as $\infty\perp q$
 and $y=o+{\fk y}+{({\fk y,y})\over 2}\,\infty\perp q$
 iff ${\fk y}\in\C$
---
we interpret this plane as the ideal boundary
 of (two) Poincar\'e half space(s) $H^3$,
 conveyed by the isometry
$$
  \{{\fk y}\in\R^3\,|\,({\fk y},q)\neq 0\} =:
   H^3 \ni {\fk y}\mapsto y :
   = {-1\over({\fk y},q)}(o+{\fk y}+{({\fk y,y})\over 2}\,\infty)
  \in K^3.
$$
Given a surface ${\fk f}:\Sigma^2\to H^3$ in (one of the)
 Poincar\'e half space(s) and its (hyperbolic) tangent plane map,
$$
  f = {-1\over({\fk f},q)}(o+{\fk f}+{({\fk f,f})\over 2}\,\infty):
   \Sigma^2\to K^3
   \enspace{\rm and}\enspace
  t:\Sigma^2\to S^{3,1},
\eqno poisurf$$
its hyperbolic Gauss map $h=q+f+t:\Sigma^2\to\lc^3$
 of (\:kilhgm:)
is a lift of a $\C$-valued map ${\fk h}:\Sigma^2\to\C$,
$$
  \langle h\rangle
   = \langle o+{\fk h}+{({\fk h,h})\over 2}\,\infty\rangle:
   \Sigma^2\to\P(\lc^3),
\eqno poihgm$$
and the enveloped horosphere congruence $s:\Sigma^2\to S^{3,1}$
 of (\:kilhgm:) can be described as in (\:eucsph:),
 in terms of
 ${\fk h}$ and
 its (Euclidean) radius function $r:\Sigma^2\to\R$:
$$
  s = t+f = h-q
    = {1\over r}\,(o+{\fk h}-r\,q+{({\fk h,h})\over 2}\,\infty):
    \Sigma^2\to S^{3,1}
\eqno poihsc$$
since $s$ has centre surface ${\fk c}={\fk h}-r\,q$.
Note that the lift $h$ of the hyperbolic Gauss map is chosen
 to be isometric, $(dh,dh)=(ds,ds)$, and
 that $(dh,dh)={1\over r^2}(d{\fk h},d{\fk h})$;
thus

\proclaim\@rconfak:Lemma 1@.
The (Euclidean) radius function $r$ of a horosphere congruence
 in Poincar\'e half space is the conformal factor relating
 the metrics
 of the ideal envelope ${\fk h}:\Sigma^2\to\C$ and
 of the sphere congruence $s:\Sigma^2\to S^{3,1}$,
$$
  (d{\fk h},d{\fk h}) = r^2\,(ds,ds).
$$

\section 2. Linear Weingarten surfaces of Bryant type

We will compute the (intrinsic) Gauss curvature
of a (Ribaucour) horosphere congruence:
it turns out that this curvature is constant exactly when
the horosphere congruence is enveloped by a Bryant type linear
Weingarten surface and (one of) its hyperbolic Gauss map(s),
 cf [\:2001jmn:,~Sect~4.3] and [\:2004gmm:, Thm 1].
This characterization of Bryant type linear Weingarten surfaces
 is the cornerstone
 of our generalization of the Bianchi-Cal\`o method.

\subsection Isothermic sphere congruences.
A surface in a space form is
 {\em linear Weingarten\/}
 if its Gauss and mean curvatures satisfy a non-trivial
 affine relation;
in particular,
a surface ${\fk f}:\Sigma^2\to H^3$ in hyperbolic space is
 {\em linear Weingarten of Bryant type with parameter $\mu$\/}
 in the sense of [\:2004gmm:]
 if there is a constant $\mu\in\R$ so that
$$
  0 = (\mu+1)\,K - 2\mu\,H + (\mu-1).
\eqno bryantlw$$

Any linear Weingarten surface comes with a distinguished pair
 of enveloped sphere congruences,
see [\:2012bjr:].
In the case (\:bryantlw:) of a Bryant type linear Weingarten
 surface (lifts of) these isothermic sphere congruences
 are given by
$$
  \sigma^+ := \tau + \varphi
   \enspace{\rm and}\enspace
  \sigma^- := {\mu+1\over 2}\,\tau + {\mu-1\over 2}\,\varphi,
\eqno isocong$$
where $\varphi=f$ and $\tau=p+t$ from (\:liemob:) and
 $f$ and $t$ denote the Killing model lift and tangent plane
 map (\:kilsurf:) of the surface in hyperbolic space.
Note that $\sigma^+=p+s$ is but a constant offset of the
 M\"obius geometric representation of the horosphere congruence of (\:kilhgm:).

The sphere congruences of (\:isocong:) can be used to
 characterize Bryant type linear Weingarten surfaces:
using Rodrigues' equations (\:rodrigues:)
 it is easily confirmed that
$$
  d\sigma^+\wedge d\sigma^-_{{\del u},{\del v}}
   = \sigma^+_u\wedge\sigma^-_v - \sigma^+_v\wedge\sigma^-_u
   = ((\mu+1)\,K-2\mu H+(\mu-1))\,\varphi_u\wedge\varphi_v,
\eqno lwcharpf$$
where $(u,v)$ denote (any) curvature line coordinates
 and $\wedge$ of $\R^{4,2}$-valued $1$-forms is the usual
 exterior product of forms with the exterior product
 on $\Lambda\R^{4,2}$ being used to multiply coefficients.

Thus the following assertion follows directly,
cf [\:2012bjr:] and [\:2026bjkp:],
or [\:2020pe3:] and [\:2025cnopry:]:

\proclaim\@lwchar:Lemma 2@.
A surface (\:kilsurf:) in hyperbolic space is linear Weingarten
 (\:bryantlw:) of Bryant type with parameter $\mu$
if and only if the sphere congruences $\sigma^\pm$ of
 (\:isocong:) satisfy
$$
  0 = d\sigma^+\wedge d\sigma^-.
$$

\subsection Horosphere congruences.
We now take a M\"obius geometric point of view and consider
 a fixed $2$-sphere $q\in S^{3,1}\subset\R^{4,1}$ and
 a congruence (\:kilhgm:) of touching spheres
$$
  s = h - q :\Sigma^2\to S^{3,1}
   \enspace\hbox{\rm with envelopes}\enspace\cases{
  \langle h\rangle:\Sigma^2\to\P(\lc\cap q^\perp), \cr
  \langle f\rangle:\Sigma^2\to\P(\lc), \cr}
\eqno mobsc$$
where the lift $f$ is adjusted as in (\:kilsurf:),
 so that $(f,q)=-1$.
To facilitate the curvature computation of the next paragraph
 we record that
 $h,f$ are isotropic parallel normal fields of $s$ as
$$
  h,f \perp ds=dh
   \enspace{\rm and}\enspace
  (h,f)=(s,f)-1\equiv -1.
\eqno paralhf$$

Next fix $\mu\in\R$ and eliminate $\tau$ from (\:isocong:)
 to define
$$
  \sigma^+ := p+s
   \enspace{\rm and}\enspace
  \sigma^- := {1+\mu\over 2}\,\sigma^+-f;
\eqno iscong$$
We observe that
 $(\sigma^+,\sigma^+)=(\sigma^+,\sigma^-)=(\sigma^-,\sigma^-)=0$
and
 $(d\sigma^+,\sigma^-)=0$,
showing that $\langle\sigma^+,\sigma^-\rangle$ defines
 a Legendre map,
 with point sphere map
$$
  (\sigma^+\wedge\sigma^-)(p)
   = (p,\sigma^+)\,\sigma^- - \sigma^+(\sigma^-,p)
   = -\sigma^-+{1+\mu\over 2}\,\sigma^+
   = f;
$$
here we use the identification
 $\Lambda^2(\R^{4,2})\cong{\fk so}(4,2)$,
 $v\wedge w\mapsto (.,v)\,w-v\,(w,.)$.
Further note that $\sigma^+\perp p-q$,
 where $p-q$ is a Lie sphere geometric incarnation
 of the ideal boundary of the hyperbolic space considered,
thus confirming that $\sigma^+$ is a horosphere congruence
 into this hyperbolic space.

\subsection Curvature computation.
We assume that $s$ immerses into $S^{3,1}$ and
 seek to determine the curvature of its induced metric.
To this end we use the decomposition
$
  \R^{4,1} = s^{(1)}\oplus(s^{(1)})^\perp
   = (\langle s\rangle\oplus ds(T\Sigma^2))
     \oplus \langle h,f\rangle
$
of $\R^{4,1}$ into the first derived bundle of $s$ and
 its normal bundle to obtain the Gauss-Weingarten equations
$$
  d = D^1 + D^\perp - (h\wedge df+f\wedge dh),
   \enspace{\rm where}\enspace
  D^1 = \nabla^s + \nabla + s\wedge ds
$$
with the radial and Levi-Civita connections
 $\nabla^s$ resp $\nabla$ of $s$.

As the radial and normal connections $\nabla^s$ and $D^\perp$
 of $s$ are flat,
with parallel vector fields $s$ resp $h$ and $f$,
the Gauss-Ricci equation reduces to
$$\matrix{
 0 &=& R^1
    + {1\over 2}[(h\wedge df+f\wedge dh)
           \wedge(h\wedge df+f\wedge dh)] \hfill\cr
   &=& R^\nabla
    + {1\over 2}[(s\wedge ds)\wedge(s\wedge ds)]
    + [(h\wedge df)\wedge(f\wedge dh)] \hfill\cr
   &=& R^\nabla
    + {1\over 2}\,ds\wedge ds
    - dh\wedge df, \hfill\cr}
\eqno sgauss$$
where the quadratic terms of the first bracket vanish
 by (\:paralhf:),
 as $h,f$ are isotropic parallel normal fields,
and the Ricci-part of the equation is trivially satisfied
 since $(dh,df)$ is symmetric.
Next we use (\:iscong:) and $dh=ds$
 to express $dh\wedge df$ in terms of $\sigma^\pm$,
$$
  dh\wedge df
   = d\sigma^+\wedge({1+\mu\over 2}\,d\sigma^+-d\sigma^-)
   = {1+\mu\over 2}\,ds\wedge ds
   - d\sigma^+\wedge d\sigma^-,
$$
which transforms the Gauss equation (\:sgauss:) further to
$
  0 = (R^\nabla - {\mu\over 2}\,ds\wedge ds)
    + d\sigma^+\wedge d\sigma^-;
$
hence
$$
  R^\nabla = {\mu\over 2}\,ds\wedge ds
   \enspace\Leftrightarrow\enspace
  0 = d\sigma^+\wedge d\sigma^-.
$$

Thus, activating \:lwchar:, we obtain a generalization
of the characterization [\:2001jmn:, Sect~4.3]
of cmc-1 surfaces in hyperbolic space
as {\em isothermic surfaces of spherical type\/},
cf [\:2004gmm:, Thm 1]:

\proclaim\@curvlw:Thm 3@.
A (regular) horosphere congruence $s:\Sigma^2\to S^{3,1}$ has
 constant (intrinsic) Gauss curvature $K_{(ds,ds)}=-\mu$
if and only if it is enveloped by
 a linear Weingarten surface $f$ of Bryant type with parameter
 $\mu$ and its hyperbolic Gauss map $h$.

Note that our computation (\:sgauss:) only uses that $h,f$
 are parallel isotropic normal fields of $s$;
hence (\:sgauss:) holds for more general Ribaucour sphere congruences
 $s$ than just horosphere congruences.

For example,
 the Ribaucour (i.e., $R^\perp=0$) sphere congruence $s$
 has flat first derived bundle $s^{(1)}$ (i.e., $R^1=0$)
if and only if
 its extended Gauss map
 $\langle f,h\rangle:\Sigma^2\to{O(4,1)\over O(3)\times O(1,1)}$
 is a {\em curved flat\/}
or, equivalently,
 cf [\:1997bjpp:] or [\:2003imdg:, \S3.3.1]:

\proclaim Cor 4.
The envelopes of a Ribaucour sphere congruence
 $s:\Sigma^2\to S^{3,1}$
are Darboux transforms of each other
if and only if $R^1=0$, that is, $K_{(ds,ds)}=1$.

\section 3. The Bianchi-Cal\`o method

We now turn to our main result:
key to the Bianchi-Cal\`o method is the interplay between
 the hyperbolic and Euclidean geometries of the Poincar\'e
 half space model of hyperbolic geometry,
both considered as subgeometries of M\"obius or Lie sphere
 geometry.
In particular, we are now prepared to
 describe the enveloped isothermic horosphere congruence $s$
 of a(ny) Bryant type linear Weingarten surface ${\fk f}$
 in terms of its holomorphic hyperbolic Gauss map ${\fk h}$ and
 the Bryant parameter $\mu$,
by providing an explicit formula for its (Euclidean)
 radius function.

From [\:2012bjr:, Sect 4] and \:curvlw:\
 we know that a linear Weingarten surface in hyperbolic space
 is of Bryant type with parameter $\mu$
 exactly when:
one of its enveloped isothermic sphere congruences $s$
 consists of horospheres,
 thus furnishing a holomorphic hyperbolic Gauss map $h$
 as its second envelope,
and has constant curvature $K_{(ds,ds)}=-\mu$,
 that is, there are (local) isometric coordinate functions
$$
  z:(\Sigma^2,(ds,ds))\to(\C,{4|dz|^2\over(1-\mu|z|^2)^2}).
\eqno isocoord$$
Adopting a Poincar\'e half space model,
 with representations (\:poihsc:) and (\:poihgm:)
 for $s$ and $h$,
the (Euclidean) radius function $r$ of $s$ occurs as
 a conformal factor, by \:rconfak::
$$
  r^2{4|dz|^2\over(1-\mu|z|^2)^2} = |d{\fk h}|^2
   = |{\fk h}'|^2|dz|^2
   \enspace\Leftrightarrow\enspace
  r = {(1-\mu|z|^2)\,|{\fk h}'|\over 2}.
\eqno radiusbc$$
Thus given a hyperbolic Gauss map ${\fk h}$,
 the horosphere congruence $s$ is known
 by (\:poihsc:)
and hence its envelopes
 $\langle h\rangle$ and $\langle f\rangle$
 are given as isotropic normal fields
 in a M\"obius geometric setting
or, equivalently, by [\:2002liro:, Prop~2.1]
 in a Euclidean setting.
In particular, no further integration is required.
Thus we obtain a Bianchi-Cal\`o method for any Bryant type
 linear Weingarten surface:

\proclaim\@bcm:Thm 5@.
Given a holomorphic map ${\fk h}:\C\supset{\cal D}\to\C$
 and $\mu\in\R$,
define a sphere congruence $s$ in terms of
 its centre surface ${\fk c}$ and
 its radius function $r$
 by
$$
  {\fk c}:=(r,{\fk h}):{\cal D}\to\R\times\C\cong\R^3
   \enspace{\sl with}\enspace
  r := {(1-\mu|z|^2)\,|{\fk h}'|\over 2};
$$
then the second envelope
 ${\fk f}:{\cal D}\to\R\times\C$
 of the (horo-)sphere congruence $s$
 is linear Weingarten of Bryant type with parameter $\mu$,
 as a surface in the Poincar\'e half space(s)
 $H^3\cong\R^\pm\times\C$.
Conversely: every linear Weingarten surface of Bryant type
 in Poincar\'e half space can locally be constructed in this way.

For cmc-$1$ surfaces ($\mu=-1$) in hyperbolic space
we recover the radius function of [\:2002liro:~(5)]
resp the centre surface
 of [\:1923bi:,~\S328; \:2002liro:~(7)]
 for the construction as described
 in [\:1923bi:,~\S496(b); \:2002liro:,~Thm~3.1].

Note that \:bcm:\ yields a method to construct a {\em unique\/}
 Bryant type linear Weingarten surface from a given
 hyperbolic Gauss map and parameter $\mu$
--- where the constructed surface depends
 on the parametrization.
For example, it is well known that Bryant type linear
 Weingarten surfaces come in parallel families,
 cf [\:2009ko:, Thms 2.3 \& 2.5];
this is realized by a reparametrization in
 the Bianchi-Cal\`o data (\:fig:~left),
$$
  \tilde\mu = \mu\,e^{-2\rho}
   \enspace{\rm and}\enspace
  \tilde{\fk h}(\tilde z) = {\fk h}(z)
   \enspace{\rm with}\enspace
  \tilde z = e^\rho z,
   \enspace{\rm hence}\enspace
  \tilde r = r\,e^{-\rho}.
$$

On the other hand, an isometric reparametrization $z=z(\tilde z)$,
 with
 ${|dz|^2\over(1-\mu|z|^2)^2}={|d\tilde z|^2\over(1-\mu|\tilde z|^2)^2}$,
does not change the radius function,
 hence not the represented surface,
 as
$$
  \tilde\mu = \mu
   \enspace{\rm and}\enspace
  \tilde{\fk h}(\tilde z) = {\fk h}(z)
   \enspace{\rm with}\enspace
  z = z(\tilde z)
   \enspace{\rm yields}\enspace
  {\tilde r(\tilde z)\over r(z)}
   = {(1-\mu|\tilde z|^2)|{\fk h}'(z)z'|\over
      (1-\mu|z|^2)|{\fk h}'(z)|}
   = 1,
$$
while a change by a more general holomorphic reparametrization,
 e.g., by a non-isometric M\"obius transformation,
yields a new surface with parameter $\mu$
 (and the same topology, see \:fig:\ right).

\vskip 5pt
\hglue .08\hsize\hbox to .84\hsize{%
\vtop{%
\pdfximage width 0.38\hsize {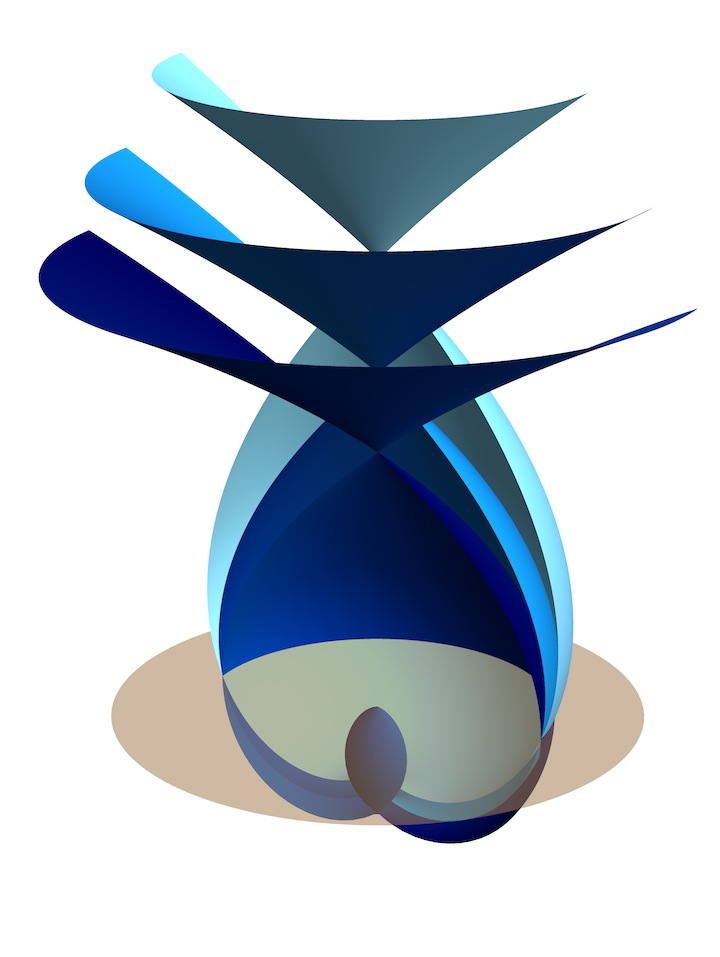}
 \pdfrefximage\pdflastximage
\pdfximage width 0.38\hsize {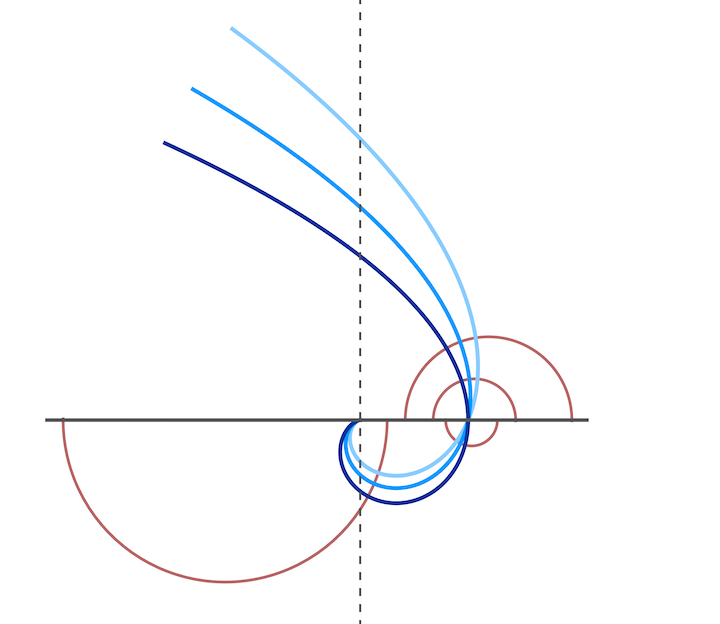}
 \pdfrefximage\pdflastximage
}%
\hfill
\vtop{\vglue 56pt%
\pdfximage width 0.46\hsize {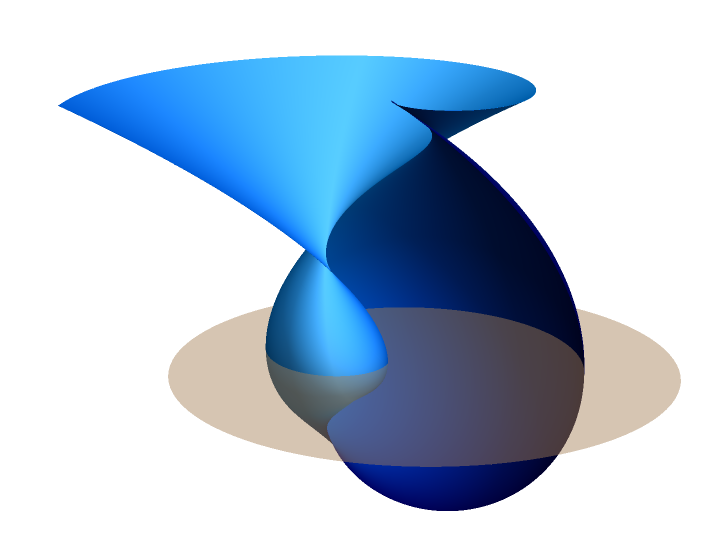}
 \pdfrefximage\pdflastximage
\pdfximage width 0.46\hsize {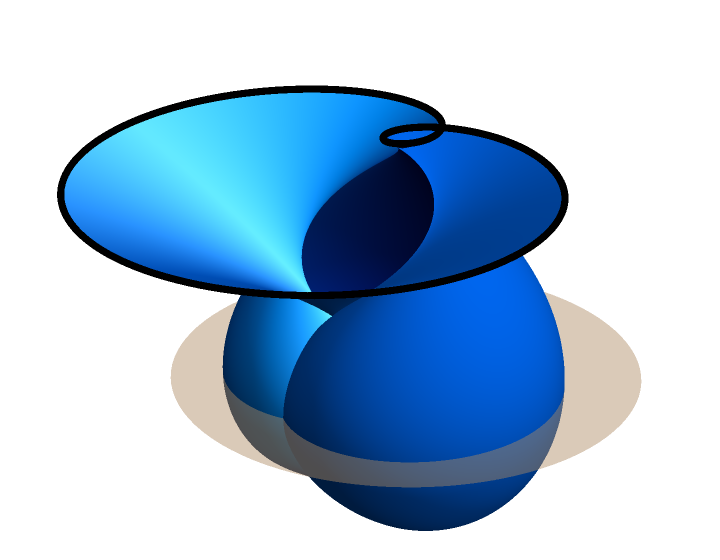}
 \pdfrefximage\pdflastximage
}%
}\vglue 0pt
\hglue .08\hsize{\vtop{\hsize=.84\hsize%
 {\bf\@fig:Fig 1@}.
  Parallel family of an HMC-1 ($\mu=1$) surface of revolution
   with hyperbolic Gauss map ${\fk h}(z)=z^2$ (top left)
   and the corresponding profile curves (bottom left);\break
  another HMC-1 surface obtained using a reparametrization
   of ${\fk h}(z)$ by a non-isometric M\"obius transformation,
   ${\fk h}(z(\tilde z))
    =\tilde{\fk h}(\tilde z)
    =({a\tilde z+b\over c\tilde z+d})^2$
   (half period top right, full period bottom right).
}}

\section References

\message{References}
\bgroup\frenchspacing\parindent=2em

\refitem 2026bjkp
 S Bentrifa, U Hertrich-Jeromin, M Kokubu, D Polly:
 {\it Ends of linear Weingarten surfaces in hyperbolic space\/};
 in preparation (2026)

\refitem 1923bi
 L Bianchi:
 {\it Lezioni di Geometria Differenziale, Vols
  \url II.1 <https://name.umdl.umich.edu/abr1998.0002.001> \&
  \url II.2 <https://name.umdl.umich.edu/abr1998.0001.001>\/};
 Nicola Zanichelli Editore, Bologna (1923)

\refitem 1929bl
 W Blaschke:
 {\it Vorlesungen \"uber Differentialgeometrie III\/};
 \url Springer Grundlehren XXIX, Berlin (1929)
  <https://doi.org/10.1007/978-3-642-50823-3>

\refitem 1997bjpp
 F Burstall, U Hertrich-Jeromin, F Pedit, U Pinkall:
 {\it Curved flats and isothermic surfaces\/};
 \url Math Z 225, 199--209 (1997)
  <https://doi.org/10.1007/PL00004308>

\refitem 2010bjr
 F Burstall, U Hertrich-Jeromin, W Rossman:
 {\it  Lie geometry of flat fronts in hyperbolic space\/};
 \url Comptes Rendus 348, 661-664 (2010)
  <https://doi.org/10.1016/j.crma.2010.04.018>

\refitem 2012bjr
 F Burstall, U Hertrich-Jeromin, W Rossman:
 {\it  Lie geometry of linear Weingarten surfaces\/};
 \url Comptes Rendus 350, 413--416 (2012)
  <https://doi.org/10.1016/j.crma.2012.03.018>

\refitem 2019bjlm
 F Burstall, U Hertrich-Jeromin, M Lara Miro:
 {\it Ribaucour coordinates\/};
 \url Beitr Alg Geom 60, 39--55 (2019)
  <https://doi.org/10.1007/s13366-018-0391-9>

\refitem 1900ca
 B Cal\`o:
 {\it Risoluzione di alcuni problemi sull'applicabilit\`a
  delle superficie\/};
 \url Ann Mat III 4, 123--130 (1900)
  <https://doi.org/10.1007/BF02419315>

\refitem 2008ce
 T Cecil:
 {\it Lie sphere geometry\/};
 \url Springer Universitext, New York (2008)
  <https://doi.org/10.1007/978-0-387-74656-2>

\refitem 2025cnopry
 J Cho, K Naokawa, Y Ogata, M Pember, W Rossman, M Yasumoto:
 {\it Discrete isothermic surfaces in Lie sphere geometry\/};
 \url Springer LNM 2375, Cham Switzerland (2025)
  <https://doi.org/10.1007/978-3-031-95592-1>

\refitem 2004gmm
 J G\'alvez, A Mart\'\i{}nez, F Mil\'an:
 {\it Complete linear Weingarten surfaces of Bryant type.
  A Plateau problem at infinity\/};
 \url Trans Amer Math Soc 356, 3405--3428 (2004)
  <https://doi.org/10.1090/S0002-9947-04-03592-5>

\refitem 2001jmn
 U Hertrich-Jeromin, E Musso, L Nicolodi:
 {\it M\"obius geometry of surfaces of constant
  mean curvature 1 in hyperbolic space\/};
 \url Ann Global Anal Geom 19, 185--205 (2001)
  <https://doi.org/10.1023/A:1010738712475>

\refitem 2003imdg
 U Hertrich-Jeromin:
 {\it Introduction to M\"obius differential geometry\/};
 \url London Math Soc Lect Note Series 300,
  Cambridge Univ Press, Cambridge (2003)
  <https://www.cambridge.org/gb/knowledge/isbn/item1150754/>

\refitem 2009ko
 M Kokubu:
 {\it Surfaces and fronts with harmonic-mean curvature one
  in hyperbolic three-space\/};
 \url Tokyo J Math 32, 177--200 (2009)
  <https://doi.org/10.3836/tjm/1249648416>

\refitem 2002liro
 L de Lima, P Roitman:
 {\it Constant mean curvature one surfaces in hyperbolic
  $3$-space using the Bianchi-Cal\`o method\/};
 \url Ann Braz Acad Sci 74, 19--24 (2002)
  <https://doi.org/10.1590/S0001-37652002000100002>

\refitem 2020pe3
 M Pember:
 {\it Lie applicable surfaces\/};
 \url Commun Anal Geom 28, 1407--1450 (2020)
  <https://doi.org/10.4310/CAG.2020.v28.n6.a5>

\refitem 2020pe
 M Pember:
 {\it Weierstrass-type representations\/};
 \url Geom Dedicata 204, 299--309 (2020)
  <https://doi.org/10.1007/s10711-019-00456-y>

\refitem 1994sm
 A Small:
 {\it Surfaces of constant mean curvature one in $H^3$ and
  algebraic curves on a quadric\/};
 \url Proc AMS 122, 1211-1220 (1994)
  <https://doi.org/10.1090/S0002-9939-1994-1209429-2>

\refitem 2020ze
 J Zeilinger:
  {\it Surfaces of constant mean curvature $1$ in hyperbolic
   $3$-space: The Bianchi-Cal\`o method and Darboux transforms\/};
  \url MSc thesis, TU Wien (2020)
   <http://hdl.handle.net/20.500.12708/79938>

\egroup

\vskip3em
\bgroup\fn[cmr7]\baselineskip=8pt
\def\addwd{\hsize=.31\hsize}
\def\udo{\vtop{\addwd
 U Hertrich-Jeromin\\
 Vienna University of Technology\\
 Wiedner Hauptstra\ss{}e 8--10/104\\
 A-1040 Vienna (Austria)\\
 udo.hertrich-jeromin@tuwien.ac.at
 }}
\def\fran{\vtop{\addwd
 F Burstall\\
 Dept of Mathematical Sciences\\
 University of Bath\\
 Bath, BA2 7AY (UK)\\
 feb@maths.bath.ac.uk
 }}
\def\gudrun{\vtop{\addwd
 G Szewieczek\\
 University Innsbruck\\
 Technikerstr 13\\
 A-6020 Innsbruck (Austria)\\
 gudrun.szewieczek@uibk.ac.at
 }}
\hbox to \hsize{\hfil\fran\hfil\udo\hfil\gudrun\hfil}\vskip 3ex
\egroup
\bye